\RequirePackage{ifpdf}
\ifpdf 
\documentclass[pdftex]{sigma}
\else
\documentclass{sigma}
\fi

\def\M{\hbox{${\cal M}$}}
\def\Ez{\hbox{$E^0$}}
\def\GR{\hbox{$\Gamma_{\mathbb{R}}$}}
\def\Gn{\hbox{$\Gamma_{n}$}}
\def\GRn{\hbox{$\Gamma_{\mathbb{R}^n}$}}

\begin{document}


\renewcommand{\PaperNumber}{020}

\FirstPageHeading

\renewcommand{\thefootnote}{$\star$}

\ShortArticleName{An Inf\/inite Dimensional Approach to the Third Fundamental Theorem of Lie}

\ArticleName{An Inf\/inite Dimensional Approach\\ to the Third Fundamental Theorem of Lie\footnote{This
paper is a contribution to the Proceedings of the Seventh
International Conference ``Symmetry in Nonlinear Mathematical
Physics'' (June 24--30, 2007, Kyiv, Ukraine). The full collection
is available at
\href{http://www.emis.de/journals/SIGMA/symmetry2007.html}{http://www.emis.de/journals/SIGMA/symmetry2007.html}}}

\Author{Richard D. BOURGIN and Thierry P. ROBART}
\AuthorNameForHeading{R.D. Bourgin and T.P. Robart}

\Address{Department of Mathematics, Howard University, Washington DC 20059, USA}
\Email{\href{mailto:rbourgin@howard.edu}{rbourgin@howard.edu}, \href{mailto:trobart@howard.edu}{trobart@howard.edu}}

\ArticleDates{Received November 02, 2007, in f\/inal form January 16, 2008; Published online February 16, 2008}

\Abstract{We revisit the third fundamental theorem of Lie (Lie III) for f\/inite dimensional Lie algebras in the context of inf\/inite dimensional matrices.}

\Keywords{Lie algebra; Ado theorem; integration; Lie group; inf\/inite dimensional matrix; representation}

\Classification{15A09; 15A29; 17A99; 17B66; 17D99; 22A22; 22A25; 22E05; 22E15; 22E45; 58H05}

\section{Introduction}

At f\/irst blush the idea of revisiting a theorem from classical (i.e.\ f\/inite dimensional) Lie group theory, no matter how pivotal, by pitching it in an
inf\/inite dimensional framework,
may seem of questionable value. However, as we shall see, the cost of this more complicated context is mitigated by the restitution of absent,
but long desired canonicity.
 This inf\/inite
dimensional setting,~${\cal M}(\infty)$, the space of inf\/inite dimensional matrices, although the analog of the space~${\cal M}(n)$ of square matrices,
provides a much richer and f\/lexible environment
in which to work, albeit with additional subtleties and complexities.  The embedding construction we will present for formal power series leads us to the
appropriate concept of invertibility in ${\cal M}(\infty)$.

Recall that the third fundamental problem of Lie is today understood as being the integrability question of a given Lie algebra into a dif\/ferential group or local group. Ado's celebrated theorem \cite{Ado 1,Ado 2} 
asserts that any f\/inite dimensional Lie algebra admits a faithful f\/inite dimensional matrix
representation. This
result provides an easy algebraic procedure for integrating any Lie algebra into a (matrix) Lie group without needing the Frobenius theorem. However, classical
matrix representations turn out to be quite limited.
Indeed, there are f\/inite dimensional Lie groups with no faithful f\/inite dimensional matrix representations (see \cite{BN} for a~modern characterization of matrix Lie groups).

This f\/inite dimensional limitation may be circumvented by introducing inf\/inite dimensional matrices. To reconstruct a given Lie group we can start from its set of right invariant vector f\/ields. The inf\/inite setting allows to encapsulate not only its Lie algebra structure but also its underlying f\/iner analytic structure; it of\/fers a way to restore in the large the above mentioned easy integration procedure.

After recalling Cartan's and Olver's versions of Lie's third fundamental theorem (Lie III) (Section~\ref{partA})
 we describe
how to embed the basic tools of classical dif\/ferential geometry into matrix form (Section~\ref{partB.1}) and of\/fer a brief survey of the main def\/initions and properties of~${\cal M}(\infty)$ (Section~\ref{partB.2}).  Finally,
we illustrate (Lie III) in this new setting (Section~\ref{partC.1}).   We get, as consequence, that all f\/inite dimensional Lie
groups, as well as local Lie groups, admit an inf\/inite
dimensional matrix representation. Technical details will be published in a forthcoming paper~\cite{BR 2}.

\section{The third fundamental theorem}\label{partA}

Let us f\/irst recall the well-known global version of Lie's third fundamental theorem~\cite{Cartan}.

\begin{theorem}[Lie III -- \'Elie Cartan]\label{theorem LieIII-Cartan}
Each Lie algebra $\cal L$ of finite dimension is the Lie algebra of a unique (up to isomorphism) connected
and simply connected Lie group $\tilde{G}$.
Moreover any connected Lie group $G$ having $\cal L$ as Lie algebra is a quotient of that simply connected Lie group~$\tilde{G}$ by a discrete
normal subgroup of its center.
\end{theorem}

We now recall a less known but still fundamental version of (Lie III). This quite recent version is due to Peter Olver~\cite{Olver 2}. It concerns the
integrability question of local Lie groups. Not all local Lie groups can be embedded into a (global) Lie group! Before stating its main theorem we follow~\cite{Olver 2} in illustrating his construction through a basic example.

 Let $G = \mathbb{R}^2$ denote the additive 2-dimensional real Lie group and let $S=\{(-1,0)\} \subset G$. Denote by $M=G  \backslash  S$ the corresponding
punctured plane. Observe that the original global action of $G$ on itself (by translation) naturally gives rise to a local action of $G$ on $M$. Finally,
construct a simply connected smooth covering space $L$ of $M$ with projection map $\pi  : L \longrightarrow M$.  Since $\pi$ is a local isomorphism, we may
lift the local action of $G$ on $M$ to $L$.  This allows us to regard the space $L$ as a topological space of mostly invertible elements containing an
identity element. In fact all but the pre-image by $\pi$ of $S^{-1}=-S=\{(1,0)\}$ are invertible elements with inverses in~$L$.

More precisely, let us start with the following def\/inition \cite[Def\/inition~2]{Olver 2} of a local Lie group.

\begin{definition}[Local Lie group]
A smooth manifold $L$ is called a {\bf local Lie group} if there exists
\begin{enumerate}
\itemsep=0pt

\item[a)] a distinguished element $e\in L$, the identity element,

\item[b)] a smooth product map $\mu: {\cal U}\rightarrow L$ def\/ined on an open subset $\cal U$ satisfying \[
    (\{e\}\times L)\cup (L\times\{e\})\subset {\cal U}\subset L\times L,
    \]
\item[c)] a smooth inversion map $\iota:{\cal V}\rightarrow L$ def\/ined on an open subset $\cal V$, $e\in{\cal V}\subset L$ such that ${\cal V}\times\iota({\cal V})\subset {\cal U}$, and $\iota({\cal V})\times{\cal V}\subset{\cal U}$, all satisfying the following properties:
     \begin{enumerate}\itemsep=0pt
     \item[(i)] {\rm identity}: $\mu(e,x)=x=\mu(x,e)$ for all $x\in L$,
     \item[(ii)] {\rm inverse}: $\mu(\iota(x),x)=e=\mu(x,\iota(x))$ for all $x\in\cal V$,
     \item[(iii)] {\rm associativity}: if $(x,y)$, $(y,z)$, $(\mu(x,y),z)$ and $(x,\mu(y,z))$ all belong to $\cal U$, then \linebreak $\mu(x,\mu(y,z))$ equals $\mu(\mu(x,y),z).$
     \end{enumerate}
\end{enumerate}
\end{definition}

With this def\/inition in mind it is possible -- see \cite{Olver 2} --  to select two open sets $\cal U$ and $\cal V$ of the covering space $L$ of $M$ to endow
$L$ with the structure of a local Lie group. Interestingly, $L$, like most local Lie groups, cannot be enlarged into a global Lie group.
The crucial observation to be made is that {\it local associativity} (assumption (iii)) may not be extended to global associativity. To see this, consider
the points  ${\bf a} =(-2,1)$, ${\bf b} =(0,-2)$ and ${\bf c} =(2,1)$ in $M$. Start at the origin~$\bf 0$ of $M$ acting f\/irst through translation by ${t\bf a}$
for $t\in[0,1]$ then, starting at $\bf a$,  acting through translation by ${t\bf b}$ for $t\in[0,1]$, and f\/inally, starting at ${\bf a}+{\bf b}$, acting
through  translation by  ${t\bf c}$ for $t\in[0,1]$. Of course we end up back at the origin since ${\bf a}+{\bf b}+{\bf c}={\bf 0}$. The triangle~$\Delta$
we have just described contains the singular point $(-1,0)$. Let us mimic this procedure on~$L$ with the corresponding lifted actions, starting at the
identity element of~$L$, $\bf \bar{0}$, in the  0$^{\rm th}$ sheet of $L$.  Denote by  ${\bf \bar{a}}$, ${\bf\bar{b}}$, ${\bf\bar{c}}$
the lifted elements on the 0$^{\rm th}$ sheet of $L$ of ${\bf a}$, ${\bf b}$ and ${\bf c}$. Then
${\bf{\bar{a}}}+({\bf{\bar{b}}}+{\bf{\bar{c}}})={\bf{\bar{0}}}$ as may be checked directly. In particular, this sum lies in the
0$^{\rm th}$ sheet of $L$.
But   $({\bf{\bar{a}}}+{\bf{\bar{b}}})+{\bf{\bar{c}}}\neq{\bf{\bar{a}}}+({\bf{\bar{b}}}+{\bf{\bar{c}}})$ since the left side lies in the {\it first} sheet
of $L$, not the 0$^{\rm th}$. In other words, global associativity has been broken!

It was P.A.~Smith and A.I.~Mal'cev who f\/irst pointed out in the 1930's \cite{Smith 1, Smith 2, Malcev} the crucial connection between globalizability and associativity.
It took sixty more years to obtain a general method for constructing local Lie groups which are not contained in any global Lie group~\cite{Olver 2}. The
local Lie group $L$ constructed above of\/fers a nice illustration of Olver's general construction:

\begin{theorem}[Lie III -- Peter Olver]\label{theorem Lie-Olver}
Let $G$ be a connected Lie group and $S\subset G$  a
closed subset of $G$ not containing the identity. Any covering manifold $\pi : L \longrightarrow G  \backslash  S$ can be regarded as a local Lie group.
Moreover each local
Lie group is contained in one of these examples.
\end{theorem}

\section[Infinite dimensional matrices]{Inf\/inite dimensional matrices}\label{partB}

There is a plethora of equivalent formats for the set $\M(\infty)$ of inf\/inite (countable) dimensional matrices. Among them the format $\mathbb{N}\times\mathbb{N}$ (with $\mathbb{N}=\left\{1,2,\ldots\right\}$ the set of positive integers and starting from the upper-left corner) distinguishes itself for its simplicity. This is the format we shall focus on throughout this paper.

\subsection{Embedding tools}\label{partB.1}

\subsubsection[The fundamental tools of differential geometry]{The fundamental tools of dif\/ferential geometry}

 Dif\/ferential geometry makes a crucial use of the Lie pseudo-groups of smooth dif\/feomorphisms in $n$ variables. Their inf\/inite jets give rise to the
groupoids \Gn\ of formal transformations (power series) in $n$-variables. For simplicity let us focus here on the one dimensional case. The general case will be presented
in \cite{BR 2}.

Let \GR\ be the Lie groupoid of formal transformations of the real line \cite{Wein}. Thus
\[
\GR = \left\{ g=\sum_{k=0}^{+\infty}a_k(x-p)^k\ |\ a_k,\, p\in \mathbb{R},\ a_1 \neq 0\right\}
\]
endowed with the double bundle
\[
\mathbb{R} \buildrel s\over\longleftarrow \GR \buildrel t\over\longrightarrow \mathbb{R}
\]
for which the source map $s$ is def\/ined by $s(g)=p$, and the target map $t$ by $t(g)=a_0$. Observe that $x$ in $g=\sum\limits_{k=0}^{+\infty}a_k(x-p)^k$ is a dummy variable which does not really need to be mentioned. It will be useful at times to make such a choice explicit. We will then write $g(x)$ instead of just~$g$ to emphasize that choice.

Suppose $g_i$, $i=1,2$ belong to \GR. Their product $g_1\cdot g_2$ is def\/ined only when $t(g_2) = s(g_1)$, and
in this case $g_1\cdot g_2 = g_1\circ g_2$,
the formal composition of series. Note that the restriction $a_1 \ne 0$
in the def\/inition of $g \in \GR$ guarantees that the formal
B\"{u}rmann--Lagrange series (see \cite{Hen} for example) for
$g^{-1}$ exists and is in \GR. Of course, for $g=\sum_{k\ge 0} a_k(x-p)^k \in$ \GR, we have $s(g^{-1})= t(g)=a_0$ and $t(g^{-1})= s(g)=p$. Evidently, it is
not possible in general
to express $g$ (or $g^{-1}$) as a power series about any point other than its source. Hence, in particular, we usually have to move the center of coordinates to go from
the power
series expansion of $g$ to that of  $g^{-1}$ when $p\ne a_0$.
This elementary observation makes it clear why, in
the representation of \GR\ in the space of
inf\/inite matrices (see below), the inverse of a given matrix may not be expressible in the same system of coordinates as is the matrix itself.  Indeed, a change of
basis is required most of the time!

\subsubsection[Matrix representation of $s^{-1}(0)$]{Matrix representation of $\boldsymbol{s^{-1}(0)}$}

There are inf\/initely many ways to represent in matrix form $s^{-1}(0)$, i.e.\ the formal (power series) transformations based at the origin of the system of coordinates.
The procedure we follow is very natural.

 For $g \in s^{-1}(0)\subset$ \GR\ and for each positive integer $m$, let $g(x)^m$ denote the $m$-fold product (as opposed to composition) of $g=g(x)$
with itself and let $u_{g(x)}$ be the column vector
\[
u_{g(x)} = (1\ \ g(x)\ \ \ g(x)^2\ \ g(x)^3\ \ldots)^T,
  \]
where $^T$ denotes transpose. Note that, with $g$ the identity transformation, we get the column vector
\[
u_x = (1 \ \ x\ \ x^2\ \ x^3\ \ldots)^T.
\]
 Finally, denote by $M_g$ the coef\/f\/icient matrix of the power series
representations of the entries in~$u_{g(x)}$ centered at $0$. Thus
\[
M_gu_x=u_{g(x)}.
\]
Then $g \buildrel\rho\over\mapsto M_g$ is a faithful algebraic representation of $s^{-1}(0)$ in
${\cal M}(\infty)$. For $g_1,g_2\in s^{-1}(0)$, their composition $g_1\circ g_2$ (def\/ined when $t(g_2)=0$) corresponds to matrix multiplication.
Moreover, if  $t(g_1)=0$ and $t(g_2)=0$ then both $M_{g_1}$ and $M_{g_2}$ are upper triangular. Otherwise said, through the inf\/inite dimensional linear
representation $\rho$, the isotropy group of the transitive Lie groupoid \GR\ is a subgroup of the nice (formal) group of upper triangular matrices \cite{Boch}.

 Consider now the element  $g = 1 - x + x^2 + \cdots + (-1)^{n}x^n + \cdots$ of\  \GR.  Note that $g \in s^{-1}(0)$, and that its associated matrix is
\begin{gather*}
M_g  = \left(
\begin{array}{cccccccc}
 1 &  0  &  0   &  0   &  0  &  0    &  0    & \dots \\
 1 &  -1 &  1   &  -1  &  1  &  -1   &  1    & \dots \\
 1 &  -2 &  3   &  -4  &  5  &  -6   &  7    &       \\
 1 &  -3 &  6   &  -10 &  15 &  -21  &  28   &       \\
 1 &  -4 &  10  &  -20 &  35 &  -56  &  84  &       \\
 1 &  -5 &  15  &  -35 &  70 &  -126 &  210 &       \\
\vdots & \vdots & & & & & & \ddots
\end{array}
\right). 
\end{gather*}
Since $g = \tau \circ h$ where $\tau(x) = 1+x$ and $h(x) = -x+x^2-x^3 + \cdots$, $M_g$ decomposes as
$M_g = M_{\tau}  M_h$ with
\begin{gather*}
M_{\tau} = \left(
\begin{array}{@{}cccccccc@{}}
 1 &  0  &   0   &  0  &  0  &  0   &  0    & \dots \\
 1 &  1  &   0   &  0  &  0  &  0   &  0    & \dots \\
 1 &  2  &   1   &  0  &  0  &  0   &  0    &       \\
 1 &  3  &   3   &  1  &  0  &  0   &  0    &       \\
 1 &  4  &   6   &  4  &  1  &  0   &  0    &       \\
 1 &  5  &   10  &  10 &  5 &   1   &  0    &       \\
\vdots & \vdots & & & & & & \ddots
\end{array}
\right) \quad \hbox{and} \quad
M_h = \left(
\begin{array}{@{}cccccccc@{}}
 1 &   0 &  0   &   0  &  0  &  0    &  0    & \dots \\
 0 &  -1 &  1   &  -1  &  1  &  -1   &  1    & \dots \\
 0 &   0 &  1   &  -2  &  3  &  -4   &  5    &       \\
 0 &   0 &  0   &  -1  &  3  &  -6   &  10   &       \\
 0 &   0 &  0   &   0  &  1  &  -4   &  10   &       \\
 0 &   0 &  0   &   0  &  0  &  -1   &  5    &       \\
\vdots & \vdots & & & & & & \ddots
\end{array}
\right)\! .
\end{gather*}
Each element of \GR\ is invertible as a groupoid element.  In particular, $g$ has an inverse in \GR\ given explicitly by
$h^{-1} \circ \tau^{-1}$.
While both $(M_{\tau})^{-1} = M_{\tau^{-1}}$ and $(M_{h})^{-1} = M_{h^{-1}} = M_h$ exist individually as matrices, we
observe that their product
$(M_h)^{-1} (M_{\tau})^{-1}$ does not make sense (i.e.\ some inf\/inite series in the matrix product do not converge). In other words,
the geometric candidate for $(M_g)^{-1}$ fails to be ``visible" (i.e.\ convergent)
 in the present chart. We are led to enlarge the classical notion of invertibility.

\subsubsection[The associated  groupoids - their matrix form]{The associated  groupoids -- their matrix form}

In the previous paragraphs we have discussed  a matrix representation for the elements of the source f\/iber $s^{-1}(0)$ of \GR. Let us extend this
representation to the entire groupoid \GR. Observe that if $g\in\GR$\ then $g = \tau_1 \circ \gamma \circ \tau_2$ where $\gamma \in s^{-1}(0)\cap t^{-1}(0)$
(i.e.\ $\gamma$ is in the isotropy group of \GR) and $\tau_1$, $\tau_2$ are translations given explicitly by $\tau_1(x) = t(g) + x$ and $\tau_2(x) = -s(g) + x$.
As pointed out above (in the special case of the translation $x \longrightarrow 1+x$), both $\tau_1$
and $\tau_2$ have matrix representations as unipotent lower triangular matrices, and $\gamma$ as an upper triangular matrix.  Hence each element of
\GR\ has a matrix representation as a formal product $L_1 U\star L_2$ in which $L_i$, $i=1,2$ are lower triangular and $U$ is upper triangular. The star
here represents a latent (undone) matrix product.

\begin{remark}
 The above mentioned latent product is not always convergent if carried out. When it is we lose crucial information by performing it: the knowledge of the source of the corresponding groupoid element. Doing so would be equivalent to prolonging an analytic germ in a~neighborhood of its point of def\/inition and then extracting the germ of the constructed section at another point. On the contrary this remark does not hold for the other product $L_1 U$. In other terms, $L_1 U$ is equivalent to $L_1\star U$. This is in itself quite surprising!
 \end{remark}

\subsection{Invertible matrices}\label{partB.2}

We just observed that the classical concept of invertibility (existence of an inverse) needs to be enlarged in order to accept in the inf\/inite
dimensional context most geometrical constructions. In this section we describe the appropriate concept together with some important properties.
See \cite{BR 1} for more details.

\subsubsection{Which matrices are really invertible?}

In most of what follows the choice $\mathbb{R}$ or $\mathbb{C}$ for scalar f\/ield is immaterial, and when it is, it is denoted $\mathbb{K}$.
As already described,
by an inf\/inite matrix we shall mean an array $A= (a^i_{\ j})$ in which the indices $i$ and $j$ run over the positive integers, $\mathbb{N}$, the
upper index denoting the row and the lower, the column. Let \Ez\ be the vector space of column vectors with f\/initely many nonzero
entries indexed by $\mathbb{N}$ whose entries are in $\mathbb{K}$. Thus \Ez\ is the inductive limit
of the ${\mathbb{K}}^n$ with the natural injections.

\begin{definition}[$\boldsymbol{\gamma}$-invertibility]
A matrix $A=(a^i_{\ j})$ is {\it $\gamma$-invertible} if both
$A$ and $A^T$ have trivial kernel as linear maps def\/ined on $E^0$.
\end{definition}

Thus $A$ is $\gamma$-invertible if and only if both its set of
columns and its set of rows are linearly independent (for {\it finite} linear combinations).
Matrices invertible in the classical sense -- i.e.\ that admit two-sided inverses -- are
$\gamma$-invertible. Moreover, for upper --
or lower -- triangular matrices, it is straightforward to check that the $\gamma$~-- and classical notions of
invertibility coincide, and each is equivalent to having all nonzero
diagonal entries. These are straightforward consequences of the def\/inition. In the next section we describe some
useful, not so obvious consequences.

\subsubsection[Main properties of $\gamma$-invertible matrices]{Main properties of $\boldsymbol{\gamma}$-invertible matrices}

Our goal here is to provide a glimpse into the rich structure of the space of $\gamma$-invertible matrices.  The $PLU$-decomposition theorem
provides justif\/ication for the def\/inition, and is a very useful characterization of such matrices. The concept of $\beta$-block mentioned below is a
useful notion in embedding groupoids of transformations in more than one variable, like \GRn\ or groupoids of formal generalized
symmetries of systems of partial dif\/ferential equations~\cite{Rob 1}.

 Recall that a permutation matrix can be characterized as follows. It is f\/illed with~$0$s and~$1$s and each of its columns and rows contains a unique~$1$.

\begin{theorem}[$\boldsymbol{PLU}$-decomposition]\label{theorem plu}
Each $\gamma$-invertible
 matrix $A$ admits a $PLU$
decomposition. Here $P$ denotes a permutation matrix,
$L$ a lower triangular unipotent matrix and $U$ an invertible
upper-triangular matrix. Conversely any such product $PLU$ is a
$\gamma$-invertible matrix.
\end{theorem}

 The connectivity of the set of $\gamma$-invertible matrices is in stark contrast to its f\/inite dimensional counterpart.

\begin{theorem}[Connectedness]\label{theorem connect}
The set \M$^{\star}$ of
$\ \gamma$-invertible matrices with real entries and indexed by $\mathbb{N}$ is arcwise
connected in the Tychonov topology. The same holds if the index set $\mathbb{N}$ for matrix indices is replaced by
$\mathbb{Z}$, or $\mathbb{R}$  by $\mathbb{C}$.
\end{theorem}

As in the f\/inite dimensional case we can characterize $\gamma$-invertible matrices using determinants.  Indeed, let
${\pi} = (\pi_1,\pi_2)$ denote an ordered pair of permutations of
$\mathbb{N}$. For any  $M\in \M(\infty)$ denote
by $\Delta^{\pi}_n(M)$ the determinant of the $n\times n$ submatrix of $M$
with rows
$\pi_1(1),\ldots,\pi_1(n)$
and columns $\pi_2(1),\ldots, \pi_2(n)$. Let $\beta:\mathbb{N}\hookrightarrow \mathbb{N} $ be an increasing injection.
Def\/ine the $\beta$-{\it block}
$\sigma$-{\it determinant of} $M$ {\it associated with} $\pi$ as the sequence
$\Delta^{\pi}_{\beta}(M)=\big\{\Delta^{\pi}_{\beta(n)}(M)\big\}_{n\in\mathbb{N}}$.
Say that this $\beta$-block
$\sigma$-determinant is {\it nonzero on} $M$,  when
$\Delta^{\pi}_{\beta(n)}(M)\neq 0$ for all $n$. When $\beta = {\rm id}_{\mathbb{N}}$ (${\rm id}_{\mathbb{N}}$ is the identity
map on $\mathbb{N}$) suppress the subscript and write more simply  $\Delta^{\pi}(M)=\left\{\Delta^{\pi}_n(M)\right\}_{n\in\mathbb{N}}$, the
$\sigma$-{\it determinant of} $M$ {\it associated} {\it with} $\pi$.

\begin{theorem}[Characterization of $\boldsymbol{\gamma}$-invertibility by $\boldsymbol{\sigma}$-determinants]\label{theorem g-inv by s-det}
Let $M\in \M$. The following are equivalent
\begin{enumerate}\itemsep=0pt
\item[(i)] $M$ is $\gamma$-invertible;
\item[(ii)] there exists a pair $\pi$, $\beta$ such that
$\Delta^{\pi}_{\beta}(M)\neq 0$;
\item[(iii)] there exists $\pi_1$ such that, with $\pi_2 = id_{\mathbb{N}}$,
$\Delta^{\pi}(M)\neq 0$;
\item[(iv)] there exists $\pi_2$ such that, with $\pi_1=id_{\mathbb{N}}$,
$\Delta^{\pi}(M)\neq 0$.
\end{enumerate}
\end{theorem}

\section{Revisiting Lie's third theorem}\label{partC}

\subsection{An algebraic approach}\label{partC.1}

In order to reformulate (Lie III, \'Elie Cartan), recall that the quotient of a connected Lie group $G$ by a discrete
normal subgroup of its center $N$ gives rise to a covering $\pi  : G \longrightarrow G /  N$ where $\pi$ is a (locally one-to-one) group homomorphism. This puts, in a neighborhood of the identity, both versions of Lie's third fundamental theorem (see Section~\ref{partA}) on the same footing. The restriction of the covering map to an appropriate (small enough) neighborhood $\bar{U}$ of the identity of the covering (local) group always leads to an {\it isomorphism} of local Lie groups
\[
\pi  : \bar{U}\subset L {\longrightarrow} U\subset G.
\]

Recall moreover that all (local) Lie groups possess a natural compatible real analytic manifold structure. Following the line of argument in Section~\ref{partB}, we may represent any local analytic transformation in $\gamma$-invertible matrix form.  This leads to

\begin{theorem}[Matrix representation]\label{theorem matrep}
Each finite dimensional Lie group or local Lie group can be faithfully represented in \M$^{\star}$, the set of
$\gamma$-invertible matrices. In this representation, group composition is represented by matrix multiplication which always turns out to be locally absolutely convergent.
\end{theorem}

The detailed proof will be published in \cite{BR 2}. It involves in an essential way the real analytic manifold structure of (local) Lie groups, the existence of a canonical chart of the f\/irst kind for local Lie groups and the Morrey--Grauert embedding theorem of analytic manifolds into a~Euclidean space \cite{Morr, Grau}. The above mentioned covering structure admits, in our context, a~counterpart involving the adjoint action \cite{BR 2}.

\subsection{Illustration}\label{partC.2}
\noindent
Let us represent $(\mathbb{R},+)$ among the $\gamma$-invertible matrices and, by using an appropriate adjoint action, map it locally isomorphically onto
the unit circle $S^1$. For $a, z \in \mathbb{C}$ denote by $\tau_a$ translation by $a$ -- that is, $\tau_a(z)=a+z$. Denote
its matrix form by $T_a$.  Thus if
\begin{gather*}
u_z = \left(
\begin{array}{c}
1\\ z\\ z^2\\ z^3\\
\vdots\\
\end{array}
\right)\qquad\hbox{then}\qquad
 \tau_a(z) \leftrightarrow T_au_z = u_{z+a},
 \end{gather*}
 where
 \[
 T_a=
 \left(
 \begin{array}{ccccc}
 1&&&&\\ a&1&&&\\ a^2&2a&1&&\\ a^3&3a^2&3a&1&\\ \vdots &&&&\ddots\\
 \end{array}
\right).
\]
 (Note that we are now working in $\Gamma_{\mathbb{C}}=\Gamma_{\mathbb{R}}\otimes \mathbb{C}$.) The embedding schema is as follows:
\[
z \buildrel \ln\over{\longrightarrow} \tilde{z} \buildrel \tau_{iy}\over{\longrightarrow} \tilde{z}+iy\buildrel \exp\over{\longrightarrow}  ze^{iy}
\]
which  transforms, by adjoint action, $\mathbb{R}$ (identif\/ied here with $\left\{iy\ |\ y\in \mathbb{R}\right\}$) onto $S^1$. In fact, an interval of $\mathbb{R}$
centered at 0 is
transformed in this way onto a neighborhood of $1\in S^1$. Finally, the extension of this local group isomorphism to the local isomophism it induces from
$(\mathbb{R},+)$ onto $S^1$ proceeds by extension, step by step, increasing the domains on which the local group operations are def\/ined, until objects
are produced which are closed under both the group operation and inversion. The details follow.

To start, using the embedding of $\Gamma_{\mathbb{R}}$ in the space of inf\/inite dimensional matrices already discussed, let the matrix forms of the functions
$z\to \ln(1+z)$ and $z\to -1+\exp(z)$, the
isotropy part of $\exp$, be
denoted with capital letters.  Thus, within the domains of absolute convergence of the power series representations centered at 0 of these analytic maps,
\begin{gather*}
\hbox{LN}= \left(
\begin{array}{@{}cccccc@{}}
1&0& 0& 0& 0&\dots\\
&1&-1/2&1/3&-1/4&\\
&& 1& -1& 11/12&\\
&&& 1& -3/2&\\
&&&& 1&\\
&&&&&\ddots\\
\end{array}
\right) \quad\hbox{and}\quad
\hbox{EXP$_0$}=\left(
\begin{array}{@{}cccccc@{}}
1&0& 0& 0& 0&\dots\\
0&1& 1/2!& 1/3!& 1/4!&\\
0&0& 1& 1& 4/3\\
0&0& 0& 1& 3/2&\\
0&0& 0& 0& 1&\\
\vdots &&&&&\ddots\\
\end{array}
\right)\!.\!
\end{gather*}
These matrices satisfy the def\/ining equations  LN$u_z = u_{\ln(1+z)} $ and EXP$_0$$u_z=u_{-1+\exp(z)}$. Thus in matrix notation, the above  scheme may be written
\[
u_z \longrightarrow T_{-1}u_z\longrightarrow (T_{iy}\hbox{LN})(T_{-1}u_z)\longrightarrow T_1[\hbox{EXP$_0$}(T_{iy}\hbox{LN})](T_{-1}u_z).
\]
Since these matrix multiplications correspond, faithfully, to the composition of the corresponding functions, in order to understand in which regions we
will have absolute convergence, we may (and do) concentrate on the functions themselves.  Thus, for example, within the radius of absolute convergence,
\begin{gather*}
(\hbox{EXP$_0$})(T_{iy}(\hbox{LN}))
 \leftrightarrow \left\{z \longrightarrow \underbrace{iy+ \left(z-\frac{z^2}{2!}+\frac{z^3}{3!}-\frac{z^4}{4!}+\cdots\right)}_\omega +  \frac{\omega^2}{2!} +  \frac{\omega^3}{3!} +\cdots\right\} \\
\qquad{} = (e^{iy}-1)+ z(e^{iy})\\
\qquad{} + z^2\left[-{{1}\over{2}}+{{1}\over{2!}}(1-iy)+{{1}\over{3!}}\left(3iy-{{3}\over{2}}(iy)^2\right)
+{{1}\over{4!}}\left(6(iy)^2-{{(iy)^3}\over{2}}\cdot 4\right)+\cdots\right]\\
\qquad{}+ z^3\left[\cdots\right]+\cdots
= (e^{iy}-1)+z(e^{iy}).
\end{gather*}
It follows that
\[
T_1[\hbox{EXP$_0$}(T_{iy}\hbox{LN})]T_{-1} \leftrightarrow e^{iy}+(z-1)e^{iy} = ze^{iy}.
\]
That is, the adjoint action $\pi (i{\mathbb{R}})\pi^{-1}$ $\left(\hbox{where } \pi(iy) = e^{iy}\right)$ sends y to the scaling function:  $z \to ze^{iy}$. Of course, in
matrix form this function has the representation
\[
\left(
\begin{array}{cccccc}
1&0&0&0&0&\dots\\
0&e^{iy}&0&0&0&\\
0&0&e^{2iy} &0&0&\\
0&0&0&e^{3iy} &0&\\
0&0&0&0&e^{4iy} &\\
\vdots &&&&&\ddots
\end{array}
\right)
\]
which may be identif\/ied with $e^{iy}\in S^1$. These calculations presuppose that all the relevant inf\/inite series converge absolutely.  Note that LN is
applied to $T_{-1}u_z$ so the function version is
$\ln(1+(z-1))$. Hence the condition for absolute convergence is $|z-1|<1$. For each f\/ixed $y\in \mathbb{R}$, the portion of $S^1$ for which $|z-1|<1$ is the arc of
$S^1$ between polar angles $-\pi/3$ and $\pi/3$. Using this representation of $iy$  as the above diagonal matrix, it is easy to expand the domain of $y$'s for
which $iy$ may be represented, by multiplication of matrices.  After one set of multiplications we will have an association between
$\left\{iy\ |\ y\in (-2\pi/3,2\pi/3)\right\}$, and within two more expansions the circle will be completely covered by an interval of $\mathbb{R}$.

 We could similarly embed Olver's basic example (see Section~\ref{partA}) in matrix form. Indeed, it suf\/f\/ices to observe that $L$ identif\/ies with the Riemann surface associated to $z\to \ln(1+z)$.

\subsection{Perspectives}\label{partC.3}

Let $G$ denote the one parameter complex group generated by the inf\/inite matrix $X$ displayed below. Thus $G = \left\{\exp{(tX)}\ |\ t\in {\mathbb{C}}\right\}$. By way of an elementary
example we will show how an adjoint action can alter $G$ in ways not possible in the f\/inite dimensional setting. Take $X$ to be the matrix with nonzero entries only in the f\/irst column given explicitly by
\[
X=\left(
\begin{array}{rrrrr}
 0 &&&& \dots\\
-1&0&&& \\
 1&&0&& \\
-1&&&0&\\
\vdots &&&&\ddots\\
\end{array}
\right).
\]
Let
\[
u_t = \exp {(tX)} =
\left(
\begin{array}{rrrrrr}
 1&&&&&\dots\\
-t&1&&&&\\
 t &&1&&&\\
-t &&& 1 &&\\
 t &&&& 1 &\\
\vdots &&&&& \ddots\\
\end{array}
\right)
\]
and take
\[
g=\left(
\begin{array}{rrrrr}
 1&-1& 1&-1 &\dots\\
&  1&-1& 1&\dots \\
&& 1&-1&\\
&&&\ddots&\ddots\\
\end{array}
\right), \qquad\hbox{so that}\qquad g^{-1}=
\left(
\begin{array}{ccccc}
1&1&&& \\
&1&1&& \\
&&1&1& \\
&&&1&\ddots \\
&&&&\ddots\\
\end{array}
\right).
\]
Then Ad$(g){G} = \left\{g^{-1}u_t g\ |\ t\in{\mathbb{C}}\right\}$.  If $M_t$ denotes $g^{-1}u_t g$ then
by direct calculation,
\[
M_t=\left(
 \begin{array}{cccccc}
1-t &t&-t&t&-t&\dots\\
0&&&&&\\
0&&&&&\\
0&&&I&&\\
0&&&&&\\
\end{array}
\right)
\]
in which the identity matrix $I$ appears starting in the second row and second column.

Observe that
\begin{gather*}
M_{t}M_{t'}  =
\left(
\begin{array}{@{}cccccc@{}}
1-t&t&-t&t&-t&\dots\\
0&1&0&0&0&\\
0&0&1&0&0&\\
0&0&0&1&0&\\
0&0&0&0&1&\\
\vdots &&&&&\ddots\\
\end{array}
\right)
\left(
\begin{array}{@{}cccccc@{}}
1-t'&t&-t'&t'&-t'&\dots\\
0&1&0&0&0&\\
0&0&1&0&0&\\
0&0&0&1&0&\\
0&0&0&0&1&\\
\vdots &&&&&\ddots\\
\end{array}
\right)\\
\phantom{M_{t}M_{t'}}{} = \!\left(
\begin{array}{@{}cccccc@{}}
(1-t)(1-t')&(1-t)t'+t & -(1-t)t' - t & (1-t)t' + t & -(1-t)t'-t\! & \dots\\
0&1&0&0&0&\\
0&0&1&0&0&\\
0&0&0&1&0&\\
0&0&0&0&1&\\
\vdots &&&&&\ddots\\
\end{array}
\right)\!\\
\phantom{M_{t}M_{t'}}{} = M_{t+t'-tt'}.
\end{gather*}
That is, addition of complex numbers:  $(t,t')\rightarrow t+t'$ has been transformed into
\[
(t,t')\rightarrow \mu(t,t') = t+t'-tt'.
\]
Observe that $M_t$ is $\gamma$-invertible unless $t=1$.  Hence the domain of $\mu$ is no longer ${\mathbb{C}}\times{\mathbb{C}}$, but instead $({\mathbb{C}}  \backslash  \{1\})\times ({\mathbb{C}}  \backslash  \{1\})$.
It is easy to check that $\{M_t\ |\ t\in {\mathbb{C}}  \backslash  \{1\}\}$ is a group under matrix multiplication isomorphic to $({\mathbb{C}}  \backslash  \{1\},\mu)$. In particular, $\mu$ is
associative.

  It is remarkable to observe that the adjoint action Ad$(g)$ of this example has immediately altered~-- without of\/fering any convergence problem~-- the topology of the group. In this case  $\mathbb{C}$ has been transformed into ${\mathbb{C}}  \backslash \{1\}$.

\pdfbookmark[1]{References}{ref}
\LastPageEnding

\end{document}